\documentclass[leqno,10pt]{amsart}

\usepackage{amsmath,amsfonts,amsthm}
\usepackage{amssymb}
\usepackage[all]{xy}
\usepackage{enumerate}
\usepackage{mathrsfs}
\usepackage{graphicx}
\usepackage{palatino}

\theoremstyle{plain}

\theoremstyle{definition}

\newcommand\dto{\dashrightarrow}

\newcommand\lto{\longrightarrow}
\newcommand\nto{\stackrel}

\def\NN{\mathbb{N}}
\def\ZZ{\mathbb{Z}}
\def\kk{\mathbb{K}}
\def\PP{\mathbb{P}}

\def\Sc{\mathscr{S}}

\def\Zc{\mathcal{Z}}

\def\Supp{\mathrm{Supp}}

\def\aaa{\mathfrak{a}}
\def\Xc{\mathscr X}
\def\G{\textnormal{\bf G}}
\def\Z.{\mathcal{Z}_\bullet}

\def\Region{{\mathfrak R _B}}
\newcommand\Spec{\textnormal{Spec}}

\newcommand\X{\textbf{X}}
\newcommand\SIR{\textnormal{Sym}_R (I)}
\newcommand\cd{\textnormal{cd}}
\def\D{\textnormal{\bf C}\textnormal{l}(\Xc)}
\newcommand\T{\textbf{T}}
\def\SSup{\mathfrak S}

\def\deg{\mathrm{deg}}

\def\im{\mathrm{im}}

\def\dim{\mathrm{dim}}

\def\mm{\mathfrak{m}}

\title[An algorithm for comp. impl. eq. of bigraded rat. surfaces]{An algorithm for computing implicit equations of bigraded rational surfaces}

\author{Nicol\'{a}s Botbol}
\address{Departamento de Matem\'atica\\
FCEN, Universidad de Buenos Aires, Argentina \\
\& Institut de Math\'ematiques de Jussieu \\
Universit\'e de P. et M. Curie, Paris VI, France \\
E-mail address: nbotbol@dm.uba.ar
}

\thanks{This work was partially supported by the project ECOS-Sud A06E04, UBACYT X064, CONICET PIP 5617 and ANPCyT PICT
17-20569, Argentina.}

\begin{document}

\begin{abstract}
In this article we show how to compute a matrix representation and the implicit equation by means of the method developed in \cite{Bot10}, using the computer algebra system Macaulay2 \cite{M2}. As it is probably the most interesting case from a practical point of view, we restrict our computations to parametrizations of bigraded surfaces. This implementation allows to compute small examples for the better understanding of the theory developed in \cite{Bot10}, and is a complement to the algorithm \cite{BD10}.
\end{abstract}

\maketitle


\section{Introduction and background}

The interest in computing explicit formulas for resultants and discriminants goes back to B\'ezout, Cayley, Sylvester and many others in the eighteenth and nineteenth centuries. 
The last few decades have yielded a rise of interest in the implicitization problem for geometric objects motivated by applications in computer aided geometric design and geometric modeling as can be seen in for example in \cite{Kal90,MC92,MC92alg,AGR95,SC95}. This phenomena has been emphasized in the latest years due to the increase of computing power, see for example \cite{AS01,co01,BCD03,BuJo03,BC05,BCJ06,BD07,Bot08,BDD08,Bot09,Bot10}. 

Under suitable hypotheses, resultants give the answer to many problems in elimination theory, including the implicitization of rational maps. 
In turn, both resultants and discriminants can be seen as the implicit equation of a suitable map (cf.\ \cite{DFS07}). Lately, rational maps appeared in computer-engineering contexts, mostly applied to shape modeling using computer-aided design methods for curves and surfaces. A very common approach is to write the implicit equation as the determinant of a matrix whose entries are easy to compute.
In general, the search of formulas for implicitization rational surfaces with base points is a very active area of research due to the fact that, in practical industrial design, base points show up quite frequently.  In \cite{MC92alg}, a perturbation is applied to resultants in order to obtain a nonzero multiple of the implicit equation. In \cite{BuJo03, BC05,BCJ06,BD07, BDD08,Bot09,Bot10} show how to compute the implicit equation as the determinant of the approximation complexes.

In \cite{Bot10} we present a method for computing the implicit equation of a hypersurface given as the image of a rational map $\phi: \Xc \dashrightarrow \PP^n$,  where $\Xc$ is an arithmetically Cohen-Macaulay toric variety. In \cite{BDD08} and \cite{Bot09}, the approach consisted in embedding the space $\Xc$ in a projective space, via a toric embedding. The need of the embedding comes from the necessity of a $\ZZ$-grading in the coordinate ring of $\Xc$, in order to study its regularity. We exploit the natural structure of the homogeneous coordinate ring of the toric variety where the map is defined. 

In \cite{Bot10} we introduce the ``good'' region in $\G$ where the approximation complex $\Z.$ and the symmetric algebra $\SIR$ has no $B$-torsion. Indeed, we define for $\gamma\in \G$, 
\[
 \Region(\gamma):=\bigcup_{0<k< \min\{m,\cd_B(R)\}} (\SSup_B(\gamma)-k\cdot \gamma)\subset \G,
\]
where $\SSup_B(\gamma):= \bigcup_{k\geq 0}(\Supp_\G(H^{k}_B(R))+k\cdot\gamma)$. This goes in the direction of proving the main theorem loc.\ cit., \cite[Thm. 4.4 and Rmk. 4.5]{Bot10}. Precisely it asserts that, when $\Xc$ is a ($n-1$)-dimensional non-degenerate toric variety over a field $\kk$, and $S$ its Cox ring (cf.\ \cite{Cox95}). For a rational map $\phi: \Xc \dto \PP^{n}$ defined by $n+1$ homogeneous elements of degree $\rho\in\D$. If $\dim(V(I))\leq 0$ in $\Xc$ and $V(I)$ is almost a local complete intersection off $V(B)$, we proved in Theorem  \cite[Thm. 4.4]{Bot10} that, 
\[
 \det((\Z.)_\gamma)=H^{\deg(\phi)}\cdot G \in \kk[\T],
\]
for all $\gamma\notin \Region(\rho)$, where $H$ stands for the irreducible implicit equation of the image of $\phi$, and $G$ is relatively prime polynomial in $\kk[\T]$. This result is a new formulation of that in \cite[Thm. 5.7]{BuJo03} and \cite[Thm. 10 and Cor. 12]{Bot09} in this new setting.

\medskip

In this article we show how to compute a matrix representation and the implicit equation with the method developed in \cite{Bot10}, following \cite{Bot09}, using the computer algebra system Macaulay2 \cite{M2}. As it is probably the most interesting case from a practical point of view, we restrict our computations to parametrizations of a bigraded surface given as the image of a rational map $\phi: \PP^1\times\PP^1 \dashrightarrow \PP^3$ given by $4$ homogeneous polinomials of bidegree $(e_1,e_2)\in\ZZ^2$. Thus, in this case, the $\Zc$-complex can be easyly computed, and the region $\Region(e_1,e_2)$ where it is acyclic is 
\begin{equation*}\label{Regione1e2}
 \Region(e_1,e_2)=(\Supp_\G(H^{2}_B(R))+(e_1,e_2))\cup(\Supp_\G(H^{3}_B(R))+2\cdot(e_1,e_2)).
\end{equation*}

This implementation allows to compute small examples for the better understanding of the theory, but we are not claiming that this implementation is optimized for efficiency; anyone trying to implement the method to solve computationally involved examples is well-advised to give more ample consideration to this issue.


\medskip

\section{Implementation in Macaulay 2}\label{example1Ch9}

Consider the rational map 
\begin{equation}
\begin{array}{ccc}
 \PP^1\times \PP^1 	&\nto{f}{\dto}	& \PP^3\\
(s:u)\times (t:v) 	&\mapsto	& (f_1:f_2:f_3:f_4)
\end{array}
\end{equation}
where the polynomials $f_i=f_i(s,u,t,v)$ are bihomogeneous of bidegree $(2,3)\in \ZZ^2$ given by
\begin{itemize}
 \item $f_1=s^2t^3+2sut^3+3u^2t^3+4s^2t^2v+5sut^2v+6u^2t^2v+7s^2tv^2+8sutv^2+9u^2tv^2+10s^2v^3+suv^3+2u^2v^3$,
 \item $f_2=2s^2t^3-3s^2t^2v-s^2tv^2+sut^2v+3sutv^2-3u^2t^2v+2u^2tv^2-u^2v^3$,
 \item $f_3=2s^2t^3-3s^2t^2v-2sut^3+s^2tv^2+5sut^2v-3sutv^2-3u^2t^2v+4u^2tv^2-u^2v^3$,
 \item $f_4=3s^2t^2v-2sut^3-s^2tv^2+sut^2v-3sutv^2-u^2t^2v+4u^2tv^2-u^2v^3$.
\end{itemize}

Our aim is to get the implicit equation of the hypersurface $\overline{\im(f)}$ of $\PP^3$.

First we load the package ``Maximal minors''
{\footnotesize \begin{verbatim}
i1 : load "maxminor.m2"
\end{verbatim}}

Let us start by defining the parametrization $f$ given by $(f_1,f_2,f_3,f_4)$.

{\footnotesize \begin{verbatim}
i2 : S=QQ[s,u,t,v,Degrees=>{{1,1,0},{1,1,0},{1,0,1},{1,0,1}}];
i3 : e1=2;
i4 : e2=3;

i5 : f1=1*s^2*t^3+2*s*u*t^3+3*u^2*t^3+4*s^2*t^2*v+5*s*u*t^2*v+6*u^2*t^2*v+ 
        7*s^2*t*v^2+8*s*u*t*v^2+9*u^2*t*v^2+10*s^2*v^3+1*s*u*v^3+2*u^2*v^3;
i6 : f2=2*s^2*t^3-3*s^2*t^2*v-s^2*t*v^2+s*u*t^2*v+3*s*u*t*v^2-3*u^2*t^2*v+
        2*u^2*t*v^2-u^2*v^3;
i7 : f3=2*s^2*t^3-3*s^2*t^2*v-2*s*u*t^3+s^2*t*v^2+5*s*u*t^2*v-3*s*u*t*v^2-
        3*u^2*t^2*v+4*u^2*t*v^2-u^2*v^3;
i8 : f4=3*s^2*t^2*v-2*s*u*t^3-s^2*t*v^2+s*u*t^2*v-3*s*u*t*v^2-u^2*t^2*v+
        4*u^2*t*v^2-u^2*v^3;
\end{verbatim}}

We construct the matrix associated to the polynomials and we relabel them in order to be able to automatize some procedures.

{\footnotesize \begin{verbatim}

i9 : F=matrix{{f1,f2,f3,f4}}

o9 = | s2t3+2sut3+3u2t3+4s2t2v+5sut2v+6u2t2v+7s2tv2+8sutv2+9u2tv2+10s2v3+
     --------------------------------------------------------------------
     suv3+2u2v3 2s2t3-3s2t2v+sut2v-3u2t2v-s2tv2+3sutv2+2u2tv2-u2v3
     --------------------------------------------------------------------
      2s2t3-2sut3-3s2t2v+5sut2v-3u2t2v+s2tv2-3sutv2+4u2tv2-u2v3
     --------------------------------------------------------------------
      -2sut3+3s2t2v+sut2v-u2t2v-s2tv2-3sutv2+4u2tv2-u2v3|

             1       4
o9 : Matrix S  <--- S
\end{verbatim}}

The reader can experiment with the implementation simply by changing the definition of the polynomials
and their degrees, the rest of the code being identical.

Let $k$ be a field. Assume $\Xc$ is the biprojective space $\PP^1_\kk\times \PP^3_\kk$. Take $R_1:=k[x_1,x_2]$, $R_2:=k[y_1,y_2,y_3,y_4]$, and $\G:=\ZZ^2$. Write $R:=R_1\otimes_k R_2$ and set $\deg(x_i)=(1,0)$ and $\deg(y_i)=(0,1)$ for all $i$. Set $\aaa_1:=(x_1,x_2)$, $\aaa_2:=(y_1,y_2,y_3,y_4)$ and define $B:=\aaa_1\cdot \aaa_2 \subset R$ the irrelevant ideal of $R$, and $\mm:= \aaa_1+\aaa_2\subset R$, the ideal corresponding to the origin in $\Spec(R)$.

By means of the Mayer-Vietoris long exact sequence in local cohomology, we have that:

\begin{enumerate}
 \item $H^2_B(R) \cong \omega_{R_1}^\vee\otimes_k\omega_{R_2}^\vee$,
 \item $H^3_B(R) \cong H^4_{\mm}(R)=\omega_{R}^\vee$,
 \item $H^\ell_B(R)=0$ for all $\ell\neq 2$ and $3$. 
\end{enumerate}

Thus, we get that:

\begin{enumerate}
 \item $\Supp_\G(H^2_B(R)) = -\NN\times \NN+(-2,0)\cup\NN\times -\NN+(0,-2)$.
 \item $\Supp_\G(H^3_B(R)) = -\NN\times -\NN+(-2,-2)$, .
\end{enumerate}

\begin{center}
 \includegraphics{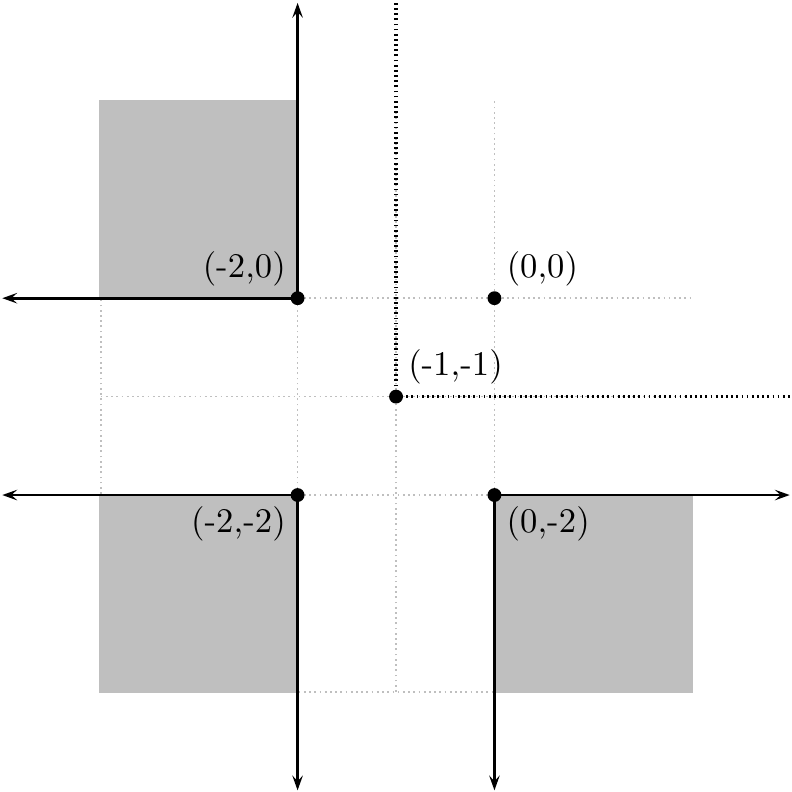}
\end{center}

And thus, 
\begin{equation*}\label{Region23}
 \Region(2,3)=(\Supp_\G(H^{2}_B(R))+(2,3))\cup(\Supp_\G(H^{3}_B(R))+2\cdot(2,3)).
\end{equation*}

Obtaining
\[
 \complement\Region(2,3)=(\NN^2+(1,5))\cup(\NN^2+(3,2)).
\]
As we can see in Example \ref{example1Ch9}, a Macaulay2 computation gives exactly this region (illustrated below) as the acyclicity region for $\Z.$.

{\footnotesize \begin{verbatim}
i10 : nu={5,3,2};
\end{verbatim}}

An alternative consists in taking 
{\footnotesize \begin{verbatim}
i10 : nu={6,1,5};
\end{verbatim}}

\begin{center}
 \includegraphics{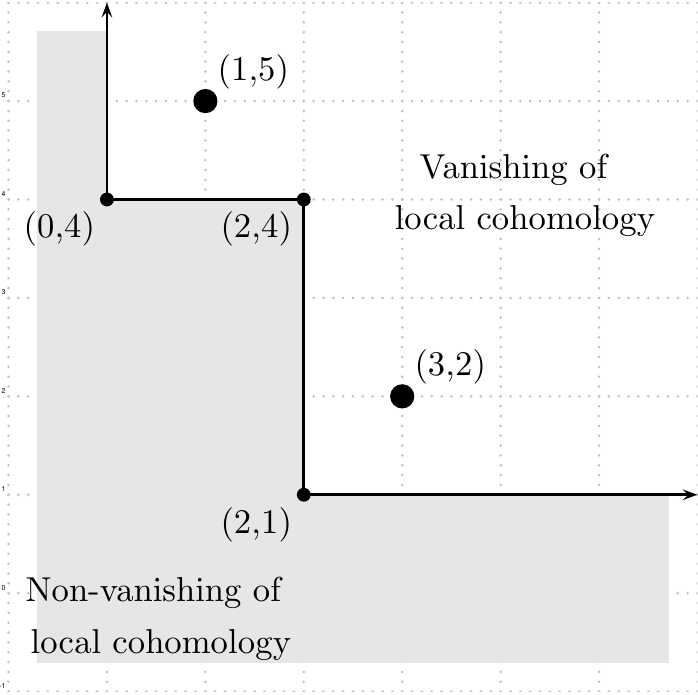}
\end{center}

Anyhow, it is interesting to test what happens in different bidegrees $\nu\in \ZZ^2$ by just replacing the desired degree in the code.
\medskip

In the following, we construct the matrix representation $M$. For simplicity, we
compute the whole module $\Zc_1$, which is not necessary as we only need the graded
part $(\Zc_1)_{\nu_0}$. In complicated examples, one should compute only this graded part
by directly solving a linear system in degree $\nu_0$. 

{\footnotesize \begin{verbatim}
i11 : Z0=S^1;
i12 : Z1=kernel koszul(1,F);
i13 : Z2=kernel koszul(2,F);
i14 : Z3=kernel koszul(3,F);

i15 : d={e1+e2,e1,e2}

i16 : hfZ0nu = hilbertFunction(nu,Z0)
o16 = 12

i17 : hfZ1nu = hilbertFunction(nu+d,Z1)
o17 = 12

i18 : hfZ2nu = hilbertFunction(nu+2*d,Z2)
o18 = 0

i19 : hfZ3nu = hilbertFunction(nu+3*d,Z3)
o19 = 0

i20 : hfnu = hfZ0nu-hfZ1nu+hfZ2nu-hfZ3nu
o20 = 0
\end{verbatim}}

Thus, when $\nu_0=(3,2)$ or $\nu_0=(1,5)$, we get a complex
\[
 (\Z.)_{\nu_0}: 0\to0\to0\to\kk[\X]^{12}\nto{M_{\nu_0}}{\lto}\kk[\X]^{12}\to0.
\]
and, hence, $\det((\Z.)_{\nu_0})=\det(M_{\nu_0})\in \kk[\X]_{12}$ is an homogeneous polynomial of degree $12$ that vanishes on the closed image of $\phi$. We compute here the degree of the MacRae's invariant which gives the degree of $\det((\Z.)_{\nu_0})$.

{\footnotesize \begin{verbatim}
i21 :hilbertFunction(nu+d,Z1)-2*hilbertFunction(nu+2*d,Z2)+
      3*hilbertFunction(nu+3*d,Z3)

o21 = 12
\end{verbatim}}

{\footnotesize \begin{verbatim}
i22 : GG=ideal F

              2 3       3   2 3   2 2        2    2 2    2   2  
o22 = ideal (s t +2s*u*t +3u t +4s t v+5s*u*t v+6u t v+7s t*v +
      ------------------------------------------------------------
              2   2   2    2 3      3   2 3    2 3   2 2       2   
      8s*u*t*v +9u t*v +10s v +s*u*v *2u v , 2s t -3s t v+s*u*t v-
      ------------------------------------------------------------
        2 2   2   2         2   2   2  2 3    2 3       3   2 2   
      3u t v-s t*v +3s*u*t*v +2u t*v -u v , 2s t -2s*u*t -3s t v+
      ------------------------------------------------------------
            2    2 2   2   2         2   2   2  2 3         3  
      5s*u*t v-3u t v+s t*v -3s*u*t*v +4u t*v -u v , -2s*u*t +
      ------------------------------------------------------------
        2 2       2   2 2   2   2         2   2   2  2 3
      3s t v+s*u*t v-u t v-s t*v -3s*u*t*v +4u t*v -u v )

o22 : Ideal of S

i23 : GGsat=saturate(GG, ideal(s,t)*ideal(u,v))

               2 2        2   2 2    2   2         2   2   2  2 3 
o23 = ideal (3s t v-3s*u*t v-u t v-3s t*v +3s*u*t*v +2u t*v -u v ,
      --------------------------------------------------------------
        2 3        2     2 2     2   2          2    2   2    2 3  
      9u t +42s*u*t v+28u t v+45s t*v -15s*u*t*v +19u t*v +30s v +
      --------------------------------------------------------------
            3    2 3       3       2   2   2         2  2   2   2 3  
      3s*u*v +13u v , s*u*t -2s*u*t v-s t*v +3s*u*t*v -u t*v , s t -
      --------------------------------------------------------------
           2    2 2    2   2         2   2   2  2 3         4  2 4 
      s*u*t v-2u t v-2s t*v +3s*u*t*v +2u t*v -u v , 30s*u*v -u v ,
      --------------------------------------------------------------
         2 4    2 4   2   3  2 4           3  2 4     2   3    2 4 
      15s v +14u v , u t*v -u v , 30s*u*t*v -u v , 15s t*v +14u v ,
      --------------------------------------------------------------
       2 2 2  2 4
      u t v -u v )

o23 : Ideal of S

i24 : degrees gens GGsat

o24 = {{{0, 0, 0}}, {{5, 2, 3}, {5, 2, 3}, {5, 2, 3}, {5, 2, 3}, {6,
      --------------------------------------------------------------
      2, 4}, {6,2, 4}, {6, 2, 4}, {6, 2, 4}, {6, 2, 4}, {6, 2, 4}}}

o24 : List

i25 : H=GGsat/GG

o25 = subquotient (| 3s2t2v-3sut2v-u2t2v-3s2tv2+3sutv2+2u2tv2-u2v3
       9u2t3+42sut2v+28u2t2v+45s2tv2-15sutv2+19u2tv2+30s2v3+3suv3+
       13u2v3 sut3-2sut2v-s2tv2+3sutv2-u2tv2 s2t3-sut2v-2u2t2v-2s2tv2+
       3sutv2+2u2tv2-u2v3 30suv4-u2v4 15s2v4+14u2v4 u2tv3-u2v4 
       30sutv3-u2v4 15s2tv3+14u2v4 u2t2v2-u2v4 |, | s2t3+2sut3+3u2t3+
       4s2t2v+5sut2v+6u2t2v+7s2tv2+8sutv2+9u2tv2+10s2v3+suv3+2u2v3 
       2s2t3-3s2t2v+sut2v-3u2t2v-s2tv2+3sutv2+2u2tv2-u2v3 2s2t3-2sut3-
       3s2t2v+5sut2v-3u2t2v+s2tv2-3sutv2+4u2tv2-u2v3 -2sut3+3s2t2v+
       sut2v-u2t2v-s2tv2-3sutv2+4u2tv2-u2v3 |)

                                1
o25 : S-module, subquotient of S

i26 : degrees gens H

o26 = {{{0, 0, 0}}, {{5, 2, 3}, {5, 2, 3}, {5, 2, 3}, {5, 2, 3}, {6,
      --------------------------------------------------------------
      2, 4}, {6,2, 4}, {6, 2, 4}, {6, 2, 4}, {6, 2, 4}, {6, 2, 4}}}

o26 : List
\end{verbatim}}

Now, we focus on the computation of the implicit equation as the determinant of the right-most map. Precisely, we will build-up this map, and later extract a maximal minor for taking its determinant. It is clear that in general it is not the determinant of the approximation complex in degree $\nu$, but a multiple of it. We could get the correct equation by taking several maximal minors and considering the gcd of the determinants. This procedure is much more expensive, hence, we avoid it.

Thus, first, we compute the right-most map of the approximation complex in degree $\nu$

{\footnotesize \begin{verbatim}
i27 : R=S[T1,T2,T3,T4];

i28 : G=sub(F,R);

              1       4
o28 : Matrix R  <--- R
\end{verbatim}} 

We compute a matrix presentation for $(\Zc_1)_\nu$ in $K_1$:
      
{\footnotesize \begin{verbatim}
i29 :Z1nu=super basis(nu+d,Z1); 

              4       12
o29 : Matrix S  <--- S

i30 : Tnu=matrix{{T1,T2,T3,T4}}*substitute(Z1nu,R);

              1       12
o30 : Matrix R  <--- R

i31 : lll=matrix {{s,t,u,v}}

o31 = | s t u v |

              1       4
o31 : Matrix S  <--- S

i32 : lll=sub(lll,R)

o32 = | s t u v |

              1       4
o32 : Matrix R  <--- R

i33 : ll={};

i34 : for i from 0 to 3 do { ll=append(ll,lll_(0,i)) }
\end{verbatim}}

Now, we compute the matrix of the map $(\Zc_1)_\nu \to A_\nu[T_1,T_2,T_3,T_4]$

{\footnotesize \begin{verbatim}
i35 : (m,M)=coefficients(Tnu,Variables=>ll,Monomials=>substitute(
            basis(nu,S),R));

i36 : M;  

              12       12
o36 : Matrix R   <--- R

i37 : T=QQ[T1,T2,T3,T4];

i38 : ListofTand0 ={T1,T2,T3,T4};

i39 : for i from 0 to 3 do { ListofTand0=append(ListofTand0,0) };

i40 : p=map(T,R,ListofTand0)

o40 = map(T,R,{T1, T2, T3, T4, 0, 0, 0, 0})

o40 : RingMap T <--- R

i41 :N=MaxCol(p(M)); 

              12       12
o41 : Matrix T   <--- T
\end{verbatim}}
 
The matrix $M$ is the desired matrix representation of the surface $\Sc$. We can continue by computing the implicit equation by taking determinant. As we mentioned, this is fairly more costly. If we take determinant what we get is a multiple of the implicit equation. One wise way for recognizing which of them is the implicit equation is substituting a few points of the surface, and verifying which vanishes.

Precisely, here there is a multiple of the implicit equation (by taking several minors we erase extra factors):
      
{\footnotesize \begin{verbatim}
i42 :Eq=det(N); factor Eq;
\end{verbatim}}

We verify the result by substituting on the computed equation, the polynomials $f_1$ to $f_4$. We verify that in this case, this is the implicit equation:
      
{\footnotesize \begin{verbatim}
i44 : use R; Eq=sub(Eq,R);
i46 : sub(Eq,{T1=>G_(0,0),T2=>G_(0,1),T3=>G_(0,2),T4=>G_(0,3)})

o46 = 0

o46 : R
\end{verbatim}}


\medskip

\subsection*{Acknowledgments\markboth{Acknowledgments}{Acknowledgments}}
 I would like to thank my two advisors: Marc Chardin and Alicia Dickenstein for the very useful discussions, ideas and suggestions.

\medskip

{\small
\def\cprime{$'$}

}

\begin{thebibliography}{AGR95}

\bibitem[AGR95]{AGR95}
Cesar Alonso, Jaime Gutierrez, and Tom\'{a}s Recio.
\newblock An implicitization algorithm with fewer variables.
\newblock {\em Comput. Aided Geom. Des.}, 12(3):251--258, 1995.

\bibitem[AS01]{AS01}
Franck Aries and Rachid Senoussi.
\newblock An implicitization algorithm for rational surfaces with no base
  points.
\newblock {\em Journal of Symbolic Computation}, 31(4):357 -- 365, 2001.

\bibitem[BC05]{BC05}
Laurent Bus{\'e} and Marc Chardin.
\newblock Implicitizing rational hypersurfaces using approximation complexes.
\newblock {\em J. Symbolic Comput}, 40(4-5):1150--1168, 2005.

\bibitem[BCD03]{BCD03}
Laurent Bus{\'e}, David Cox, and Carlos D'Andrea.
\newblock Implicitization of surfaces in {${\PP}\sp 3$} in the presence of base
  points.
\newblock {\em J. Algebra Appl}, 2(2):189--214, 2003.

\bibitem[BCJ09]{BCJ06}
Laurent Bus{\'e}, Marc Chardin, and Jean-Pierre Jouanolou.
\newblock Torsion of the symmetric algebra and implicitization.
\newblock {\em Proc. Amer. Math. Soc}, 137(6):1855--1865, 2009.

\bibitem[BD07]{BD07}
Laurent Bus{\'e} and Marc Dohm.
\newblock Implicitization of bihomogeneous parametrizations of algebraic
  surfaces via linear syzygies.
\newblock In {\em I{SSAC} 2007}, pages 69--76. ACM, New York, 2007.

\bibitem[BD10]{BD10}
Nicol{\'a}s Botbol and Marc Dohm.
\newblock A package for computing implicit equations of parametrizations from
  toric surfaces.
\newblock {\em arXiv:1001.1126}, 2010.

\bibitem[BDD09]{BDD08}
Nicol{\'a}s Botbol, Alicia Dickenstein, and Marc Dohm.
\newblock Matrix representations for toric parametrizations.
\newblock {\em Comput. Aided Geom. Design}, 26(7):757--771, 2009.

\bibitem[Bot09]{Bot08}
Nicol{\'a}s Botbol.
\newblock The implicitization problem for {$\phi\colon\mathbb {P}^n
  \dashrightarrow (\mathbb P^1)^{n+1}$}.
\newblock {\em J. Algebra}, 322(11):3878--3895, 2009.

\bibitem[Bot10a]{Bot09}
Nicol{\'a}s Botbol.
\newblock Compactifications of rational maps and the implicit equations of
  their images.
\newblock {\em To appear in J. Pure and Applied Algebra. arXiv:0910.1340},
  2010.

\bibitem[Bot10b]{Bot10}
Nicol{\'a}s Botbol.
\newblock Implicit equation of multigraded hypersurfaces.
\newblock {\em arXiv:1007.3437}, 2010.

\bibitem[BJ03]{BuJo03}
Laurent Bus{\'e} and Jean-Pierre Jouanolou.
\newblock On the closed image of a rational map and the implicitization
  problem.
\newblock {\em J. Algebra}, 265(1):312--357, 2003.

\bibitem[Cox95]{Cox95}
David~A Cox.
\newblock The homogeneous coordinate ring of a toric variety.
\newblock {\em J. Algebraic Geom}, 4(1):17--50, 1995.

\bibitem[Cox01]{co01}
David~A Cox.
\newblock Equations of parametric curves and surfaces via syzygies.
\newblock In {\em Symbolic computation: solving equations in algebra, geometry,
  and engineering ({S}outh {H}adley, {MA}, 2000)}, volume 286 of {\em Contemp.
  Math}, pages 1--20. Amer. Math. Soc, Providence, RI, 2001.

\bibitem[DFS07]{DFS07}
Alicia Dickenstein, Eva~Maria Feichtner, and Bernd Sturmfels.
\newblock Tropical discriminants.
\newblock {\em J. Amer. Math. Soc}, 20(4):1111--1133 (electronic), 2007.

\bibitem[GS]{M2}
D.~R. Grayson and M.~E. Stillman.
\newblock Macaulay 2, a software system for research in algebraic geometry.
\newblock {\em http://www.math.uiuc.edu/Macaulay2/}.

\bibitem[Kal91]{Kal90}
Michael Kalkbrenner.
\newblock Implicitization of rational parametric curves and surfaces.
\newblock In {\em AAECC-8: Proceedings of the 8th International Symposium on
  Applied Algebra, Algebraic Algorithms and Error-Correcting Codes}, pages
  249--259, London, UK, 1991. Springer-Verlag.

\bibitem[MC92a]{MC92alg}
Dinesh Manocha and John~F. Canny.
\newblock Algorithm for implicitizing rational parametric surfaces.
\newblock {\em Computer Aided Geometric Design}, 9(1):25 -- 50, 1992.

\bibitem[MC92b]{MC92}
Dinesh Manocha and John~F. Canny.
\newblock Implicit representation of rational parametric surfaces.
\newblock {\em Journal of Symbolic Computation}, 13(5):485 -- 510, 1992.

\bibitem[SC95]{SC95}
Tom Sederberg and Falai Chen.
\newblock Implicitization using moving curves and surfaces.
\newblock 303:301--308, 1995.

\end{thebibliography}
\end{document}